\magnification=1200

\def\Q{{\bf {Q}}}

\def\Z{{\bf Z}}

\def\house#1{\setbox1=\hbox{$\,#1\,$}%
\dimen1=\ht1 \advance\dimen1 by 2pt \dimen2=\dp1 \advance\dimen2 by 2pt
\setbox1=\hbox{\vrule height\dimen1 depth\dimen2\box1\vrule}%
\setbox1=\vbox{\hrule\box1}%
\advance\dimen1 by .4pt \ht1=\dimen1
\advance\dimen2 by .4pt \dp1=\dimen2 \box1\relax}

  \def\eps{{\varepsilon}}

  \def\noi{\noindent}

\def\build#1_#2^#3{\mathrel{\mathop{\kern 0pt#1}\limits_{#2}^{#3}}}

\def\date {le\ {\the\day}\ \ifcase\month\or janvier
\or fevrier\or mars\or avril\or mai\or juin\or juillet\or
ao\^ut\or septembre\or octobre\or novembre
\or d\'ecembre\fi\ {\oldstyle\the\year}}

\font\fivegoth=eufm5 \font\sevengoth=eufm7 \font\tengoth=eufm10

\newfam\gothfam \scriptscriptfont\gothfam=\fivegoth
\textfont\gothfam=\tengoth \scriptfont\gothfam=\sevengoth

\def\pro{\noindent {\bf Proof : }}

\def\smallsquare{\vbox{\hrule\hbox{\vrule height 1 ex\kern 1 ex\vrule}\hrule}}
\def\cqfd{\hfill \smallsquare\vskip 3mm}


\centerline{}

\vskip 4mm

\centerline{
\bf On the complexity of algebraic numbers II. Continued fractions}

\vskip 8mm
\centerline{Boris A{\sevenrm DAMCZEWSKI}  \
\& \ Yann B{\sevenrm UGEAUD} 
\footnote{*}{Supported by the Austrian Science
Fundation FWF, grant M822-N12. } }

\vskip 6mm

\vskip 8mm

\centerline{\bf 1. Introduction}

\vskip 6mm

Let $b \ge 2$ be an integer. \'Emile Borel [9] conjectured that 
every real irrational algebraic number $\alpha$ should satisfy
some of the laws shared by almost all real numbers with
respect to their $b$-adic expansions. 
Despite some recent progress [1, 3, 7], we are still very 
far away from establishing such a strong result.
In the present work, we are concerned with a similar question, where the
$b$-adic expansion of $\alpha$ is 
replaced  by its sequence of partial quotients.
Recall that the continued fraction expansion of an
irrational number $\alpha$ is eventually periodic if, and only if, 
$\alpha$ is a quadratic irrationality.
However, very little is known regarding the size of the partial quotients
of algebraic real numbers of degree at least three. Because of some 
numerical evidence and a belief that these numbers behave 
like most of the numbers in this respect, 
it is often conjectured that their partial quotients form an unbounded 
sequence, but we seem to be very far away from a proof (or a disproof). 
Apparently, Khintchine [16] 
was the first to consider such a question 
(see [4, 27, 29]
for surveys including a discussion on this problem). Although almost 
nothing has been proved yet in this direction, some more general 
speculations are due to Lang [17], including the fact that 
algebraic numbers of degree at least three should behave like most 
of the numbers with respect to the Gauss--Khintchine--Kuzmin--L\'evy laws.

More modestly, we may expect that 
if the sequence of partial quotients of an irrational 
number $\alpha$ is,
in some sense, `simple', then $\alpha$ is either quadratic or transcendental. 
The term `simple' can of course lead to many interpretations.  
It may denote real numbers whose continued fraction expansion 
has some regularity, or 
can be produced by a simple algorithm (by a simple Turing machine, for
example), or arises
from a simple dynamical system... 
The main results of the present work are two new combinatorial
transcendence criteria, which considerably
improve upon those from [5, 13, 8]. 
It is of a particular interest that 
such criteria naturally yield, in a unified way, several 
new results on the different approaches  
of the above mentioned notion of simplicity/complexity for the 
continued fraction expansions of algebraic real numbers of degree
at least three.

This article is organized as follows. Section 2 is devoted to the
statements of our two transcendence criteria. 
Several applications of them
are then briefly discussed in Section 3. All the proofs are postponed 
to Sections 4 and 5.

\vskip 9mm

\centerline{\bf 2. Transcendence criteria for stammering continued fractions}

\vskip 6mm

Before stating our theorems, we need to introduce some notation.
Let ${\cal A}$ be a given set, not necessarily finite. 
The length of a word
$W$ on the alphabet ${\cal A}$, that is, the number of letters
composing $W$, is denoted by $\vert W\vert$.
For any positive integer $\ell$, we write
$W^\ell$ for the word $W\ldots W$ ($\ell$ times repeated concatenation
of the word $W$). More generally, for any positive rational number
$x$, we denote by $W^x$ the word
$W^{[x]}W'$, where $W'$ is the prefix of
$W$ of length $\left\lceil(x- [x])\vert W\vert\right\rceil$. 
Here, and in all what follows, $[y]$ and
$\lceil y\rceil$ denote, respectively, the integer part and the upper
integer part of the real number $y$. 
Let ${\bf a}=(a_\ell)_{\ell \ge 1}$ 
be a sequence of elements from ${\cal A}$,
that we identify with the infinite word $a_1 a_2 \ldots a_\ell \ldots$
Let $w$ be a rational number with $w>1$.
We say that ${\bf a}$ 
satisfies Condition $(*)_w$ if ${\bf a}$ is not
eventually periodic and if there exists 
a sequence of finite words $(V_n)_{n \ge 1}$ such that:

\medskip

\item{\rm (i)} For any $n \ge 1$, the word $V_n^w$ is a prefix
of the word ${\bf a}$;

\smallskip

\item{\rm (ii)} The sequence $(\vert V_n\vert)_{n \ge 1}$ is 
increasing.

\medskip

Roughly speaking, ${\bf a}$ satisfies Condition $(*)_w$
if ${\bf a}$ is not eventually periodic and if there exist infinitely 
many `non-trivial' repetitions (the size of which is measured
by $w$) at the beginning of the infinite word $a_1 a_2 \ldots a_\ell \ldots$

Our transcendence criterion for `purely' stammering continued
fractions can be stated as follows.

\proclaim Theorem 1. 
Let ${\bf a}=(a_\ell)_{\ell \ge 1}$ be a sequence of positive integers.
Let $(p_\ell/q_\ell)_{\ell \ge 1}$ denote the sequence of convergents to 
the real number
$$
\alpha:= [0; a_1, a_2, \ldots, a_\ell, \ldots].
$$ 
If there exists a rational number $w \ge 2$ such that
${\bf a}$ satisfies Condition $(*)_w$, then $\alpha$ is transcendental.
If there exists a rational number $w>1$ such that
${\bf a}$ satisfies Condition $(*)_w$, and if
the sequence $(q_\ell^{1/\ell})_{\ell \ge 1}$ is bounded
(which is in particular the case when the sequence ${\bf a}$ is bounded),
then $\alpha$ is transcendental.

The main interest of the first statement of Theorem 1 is that there
is no condition on the growth of the sequence $(q_\ell)_{\ell \ge 1}$.
Apparently, this fact has not been observed previously.
The second statement of
Theorem 1 improves upon Theorem 4 from [5], which
requires, together with some extra rather constraining hypotheses, 
the stronger assumption $w > 3/2$. 
The condition that the sequence $(q_\ell^{1/\ell})_{\ell \ge 1}$ 
has to be bounded is in general very easy to check, and is
not very restrictive, since it is
satisfied by almost all real numbers (in the sense of the Lebesgue measure).
Apart from this assumption, Theorem 1 does not depend on the size 
of the partial quotients of $\alpha$. 
This is in a striking contrast to 
all previous results [5, 13, 8], 
in which, roughly speaking, the size $w$ of the
repetition is required to be all the
more large than the partial quotients are big.
Unlike these results, our Theorem 1 can be easily
applied even if $\alpha$ has unbounded partial quotients.

Unfortunately, in the statement of Theorem 1,
the repetitions must appear at the very beginning of ${\bf a}$.
Results from [13] allow a shift, whose length, 
however, must be controlled in terms of the size of the repetitions. 
Similar results cannot be deduced from our Theorem 1. 
However, many ideas from the proof of Theorem 1 can be used
to deal also with this situation, under some extra
assumptions, and to improve upon the transcendence criterion
from [13].

Keep the notation introduced at the beginning of this section.
Let $w$ and $w'$ be non-negative rational numbers with $w>1$.
We say that ${\bf a}$ 
satisfies Condition $(**)_{w, w'}$ if ${\bf a}$ is not
eventually periodic and if there exist 
two sequences of finite words $(U_n)_{n \ge 1}$,   
$(V_n)_{n \ge 1}$ such that:

\medskip

\item{\rm (i)} For any $n \ge 1$, the word $U_nV_n^w$ is a prefix
of the word ${\bf a}$;

\smallskip

\item{\rm (ii)} The sequence
$({\vert U_n\vert} / {\vert V_n\vert})_{n \ge 1}$ is bounded 
from above by $w'$;

\smallskip

\item{\rm (iii)} The sequence $(\vert V_n\vert)_{n \ge 1}$ is 
increasing.

\medskip

We are now ready to state our transcendence criterion for (general)
stammering continued fractions.

\proclaim Theorem 2. 
Let ${\bf a}=(a_\ell)_{\ell \ge 1}$ be a sequence of positive integers.
Let $(p_\ell/q_\ell)_{\ell \ge 1}$ denote the sequence of convergents to 
the real number
$$
\alpha:= [0; a_1, a_2, \ldots,a_n,\ldots].
$$ 
Assume that the sequence $(q_\ell^{1/\ell})_{\ell \ge 1}$ is bounded
and set $M = \limsup_{\ell \to + \infty} \, q_{\ell}^{1/\ell}$ and
$m = \liminf_{\ell \to + \infty} \, q_{\ell}^{1/\ell}$.
Let $w$ and $w'$ be non-negative real numbers with
$$
w > (2 w'+1) {\log M \over \log m} - w'. \eqno (1)
$$
If ${\bf a}$ satisfies Condition $(**)_{w, w'}$, 
then $\alpha$ is transcendental.

We display an immediate consequence of Theorem 2.

\proclaim Corollary 1.
Let ${\bf a}=(a_\ell)_{\ell \ge 1}$ be a sequence of positive integers.
Let $(p_\ell/q_\ell)_{\ell \ge 1}$ denote the sequence of convergents to 
the real number
$$
\alpha:= [0; a_1, a_2, \ldots, a_\ell ,\ldots].
$$ 
Assume that the sequence $(q_{\ell}^{1/\ell})_{\ell \ge 1}$ converges.
Let $w$ and $w'$ be non-negative real numbers with $w > w'+1$.
If ${\bf a}$ satisfies Condition $(**)_{w, w'}$, 
then $\alpha$ is transcendental.

Our Theorem 2 improves Theorem 6.3 of Davison [13].
Indeed, to apply his transcendence criterion, $w$ and $w'$ must
satisfy 
$$
w > \biggl( 2 w' + {3 \over 2} \biggr) \, {\log M \over \log m},
$$
which is a far stronger condition than (1).

Theorems 1 and 2 yield many new results 
that could not be obtained with the earlier transcendence criteria.  
Some of them are stated in Section 3, while many others will be
given in a subsequent work [2].
Theorems 1 and 2 are of the same spirit as the following
result, established in [1, 3], 
and which deals with the transcendence of $b$-adic expansions.

\proclaim Theorem ABL.
Let $b \ge 2$ be an integer.
Let ${\bf a}=(a_\ell)_{\ell \ge 1}$ be a sequence of integers
in $\{0, \ldots, b-1\}$.
Let $w$ and $w'$ be non-negative rational numbers with $w > 1$.
If ${\bf a}$ satisfies Condition $(**)_{w, w'}$, 
then the real number $\sum_{\ell \ge 1} \, a_\ell / b^\ell$ is transcendental.

Theorem ABL is as strong for `purely' stammering sequences as for
general stammering sequences, provided that the repetitions do not
occur too far away from the beginning of the infinite word. 
Unfortunately, we are
unable to replace in Theorem 2 the assumption `$w > w' + 1$' by the 
weaker one `$w > 1$', occurring in Theorem ABL.

The main tool for the proofs of Theorems 1 and 2, given in Section 4,
is the Schmidt Subspace Theorem [25, 26]. This 
(more precisely, a $p$-adic version of it) is also
the key auxiliary result for establishing Theorem ABL.

\vskip 6mm

{\bf 3. Applications to the complexity of algebraic continued fractions}

\vskip 6mm

Our transcendence criteria apply to establish that several well-known
continued fractions are transcendental, including the
Thue--Morse continued fraction (whose transcendence was first
proved by M. Queff\'elec [21]), the
Rudin--Shapiro continued fraction, folded continued
fractions, continued fractions arising from perturbed
symmetries (these sequences were introduced 
by Mend\`es France [18]),
continued fractions considered by
Davison [13] and Baxa [8], etc.
These applications are discussed in 
details in [2], where complete proofs are given. 
We only focus here on applications related to our main problem, that is, to 
the complexity of algebraic numbers 
with respect to their continued fraction expansions.

\vskip 6mm

\noindent{\bf 3.1. An algorithmic approach}

\vskip 6mm

We first briefly discuss how the complexity of the continued fraction of 
real numbers can be interpreted in an algorithmic way. Following the 
pioneering work of Turing [28], a sequence is said to be  
computable if there exists a Turing machine capable to produce  
successively its terms. Later, Hartmanis and Stearns [15]
proposed to emphasize the quantitative 
aspect of this notion, and to take into
account the number $T(n)$ of operations  
needed by a (multitape) Turing machine to produce 
the first $n$ terms of the sequence. In this regard, a real number 
is considered all the more simple 
than its continued fraction expansion can be produced very fast by 
a Turing machine. 

Finite automata are one of the most basic models of computation and take thus 
place at the bottom of the hierarchy of Turing machines. 
In particular, such machines produce 
sequences in real time, that is, with $T(n)=O(n)$. 
An infinite sequence ${\bf a}=(a_n)_{n\geq 0}$ is 
said to be generated by a $k$-automaton if $a_n$ is a 
finite-state function of the base-$k$ 
representation of $n$. 
This means that there exists a finite automaton starting 
with the $k$-ary expansion of $n$ as input and producing the term $a_n$ as 
output. A nice reference on this topic is the 
 book of Allouche and Shallit [6]. 
As a classical example of a sequence generated by a $2$-automaton, 
we mention the famous 
binary Thue--Morse sequence ${\bf a}=(a_n)_{n\geq 0}= 0110100110010 \ldots$
This sequence is defined 
as follows: $a_n$ is equal to $0$ (resp. to $1$)
if the sum of the digits in the 
binary expansion of $n$ is even (resp. is odd). 
In view of the above discussion, we may expect that finite automata are 
`too simple' Turing machines to produce the continued fraction expansion 
of algebraic numbers that are neither rationals nor quadratics. 

\proclaim Problem 1. 
Do there exist algebraic numbers of degree at least three 
whose continued fraction expansion can be produced by a finite automaton?

Thanks to Cobham [11], we know that sequences generated by finite 
automata can be characterized in terms of iterations of morphisms of free 
monoids generated by finite sets. We recall now this useful description. 
For a finite set ${\cal A}$, let ${\cal A}^*$ denote the free monoid 
generated by ${\cal A}$. The empty word is the neutral element 
of ${\cal A}^*$. Let ${\cal A}$ and ${\cal B}$ be two finite sets. An 
application from ${\cal A}$ to ${\cal B}^*$ can be uniquely extended to a 
homomorphism between the free monoids ${\cal A}^*$ and ${\cal B}^*$. 
Such a homomorphism is
called a morphism from ${\cal A}$ to ${\cal B}$. 
If there is a positive integer
$k$ such that each element of ${\cal A}$ is mapped to a word of
length $k$, then the morphism is called $k$-uniform or simply
uniform. 
Similarly, an 
application from ${\cal A}$ to ${\cal B}$ can be uniquely extended to a  
homomorphism between the free monoids ${\cal A}^*$ and ${\cal B}^*$. 
Such an application is called a coding (the term
`letter-to-letter' morphism is also used in the literature).

A morphism 
$\sigma$ from ${\cal A}$ into itself is said to 
be prolongable if there exists a 
letter $a$ such that $\sigma(a)=aW$, where the word $W$ is such 
that $\sigma^n(W)$ is a non-empty word for every $n\geq 0$. 
In that case, the 
sequence of finite words $(\sigma^n(a))_{n\geq 1}$ converges in 
${\cal A}^{\Z_{\ge 0}}$ (endowed with the product topology of the discrete 
topology on each copy of ${\cal A}$) 
to an infinite word ${\bf a}$. This infinite word is clearly 
a fixed point for $\sigma$. We say that a sequence ${\bf b}$ 
is generated by the morphism 
$\sigma$ if there exists 
a coding $\varphi$ such that ${\bf b}=\varphi({\bf a})$ . 
If, moreover, every letter appearing in ${\bf a}$
occurs at least twice, then we say that
${\bf b}$ is generated by a recurrent morphism. 
If the alphabet ${\cal A}$ has only two letters, then we say that
${\bf b}$ is generated by a binary morphism. 
Furthermore, if $\sigma$ is uniform, then we say that 
${\bf b}$ is generated by a uniform morphism. 

For instance, the Fibonacci morphism $\sigma$ defined on the alphabet 
$\{0,1\}$ by $\sigma(0)=01$ and $\sigma(1)=1$ is a binary, 
recurrent and non-uniform
morphism which generates the celebrated Fibonacci infinite 
word 
$$
{\bf a}=\lim_{n\to +\infty}\sigma^n(0)=010010100100101001\ldots
$$

Uniform morphisms and automatic sequences are strongly connected, as
shown by the following result of Cobham [11].

\proclaim Theorem  (Cobham).
A sequence can be generated by a finite automaton if, and only if, it is
generated by a uniform morphism.

This useful description gives rise to   
the following challenging question.

\proclaim Problem 2. 
Do there exist algebraic numbers of degree at least three 
whose continued fraction expansion is generated by a morphism?

Our main contribution towards both problems is the following
result.

\proclaim Theorem 3. The continued fraction expansion of an algebraic number 
of degree at least three cannot be generated by a recurrent morphism. 

The class of primitive morphisms
has been extensively studied.
In particular, Theorem 3 fully solved a question 
studied by M. Queff\'elec [22]. 
We display the following direct consequence of Theorem 3.

\proclaim Corollary 1. The continued fraction expansion of an algebraic number 
of degree at least three cannot be generated by a binary morphism. 

Indeed, it is easy to see that binary morphims are either recurrent or 
they generate only eventually periodic sequences.

\vskip 6mm

\noindent{\bf 3.2. A dynamical approach}

\vskip 6mm

In this Section, we discuss the notion of complexity of the
continued fraction expansion of a real number from a dynamical point
of view.

Let ${\cal A}$ be a given set, finite or not. 
A subshift on ${\cal A}$ is a symbolic dynamical system $(X,S)$, 
where $S$ is the
classical shift transformation defined from ${\cal A}^{\Z_{\ge 1}}$ into
itself by $S((a_n)_{n\geq 1})=(a_n)_{n\geq 2}$ and
$X$ is a subset of ${\cal A}^{\Z_{\ge 1}}$ such that $S(X)\subset X$. 
With an infinite
sequence ${\bf a}$ in ${\cal A}^{\Z_{\ge 1}}$, we associate the subshift 
${\cal X}_{\bf a}=(X,S)$, where 
$X := \overline{{\cal O}({\bf a})}$ denotes the closure of the orbit of the
sequence ${\bf a}$ under the action of $S$. The complexity function 
$p_{\bf a}$ of
a sequence ${\bf a}$ associates with any positive
integer $n$ the number $p_{\bf a}(n)$ of distinct blocks of $n$ consecutive
letters occurring in it. More generally, 
the complexity function 
$p_{\cal X}$ of a subshift ${\cal X} = (X, S)$ associates with any positive
integer $n$ the number $p_{\cal X}(n)$ of distinct blocks of $n$ consecutive
letters occurring in at least one element of $X$. 

With a subshift ${\cal X}=(X,S)$ on $\Z_{\ge 1}$ one can associate 
the set ${\cal C}_{\cal X}$ defined by
$$
{\cal C}_{\cal X}=\left\{\alpha\in (0,1), \;
\alpha=[0;a_1,a_2\ldots]\;\hbox{ such that $(a_n)_{n\geq 1}\in {\cal X}$}
\right\}.
$$
In particular, if a real number $\alpha$ lies in ${\cal C}_{\cal X}$, 
then this is also the case for any $\beta$ in
${\cal C}_{\alpha}:=\overline{(T^n(\alpha))_{n\geq 0}}$, 
where $T$ denotes the Gauss map, 
defined from $(0,1)$ into itself by
$T(x)=\{{1\over x}\}$. 
Indeed, we clearly
have $T([0;a_1,a_2,\ldots])=[0;a_2,a_3,\ldots]$. A way to 
investigate the question of the complexity 
of the continued fraction expansion
of $\alpha$ is to determine the behaviour 
of the sequence $(T^n(\alpha))_{n\geq 0}$ or, equivalently, to determine
the structure of the underlying dynamical system 
$({\cal C}_{\alpha}, T)$, 
Roughly speaking, we can consider that the larger ${\cal C}_{\alpha}$ is, 
the more complex is the continued fraction expansion of $\alpha$.

Thus, if the symbolic dynamical system ${\cal X}$ 
has a too simple structure, for
instance if it has a low complexity, we
can expect that no algebraic number of degree at least three lies in 
the set ${\cal C}_{\cal X}$.  

\proclaim Problem 3. Let ${\cal X }$ be a subshift 
on $\Z_{\ge 1}$ with sublinear
complexity, that is, whose
complexity function satisfies 
$p_{\cal X}(n) \le M n$ for some absolute constant $M$
and any positive integer $n$. Does the set ${\cal
C}_{\cal X}$ only contain quadratic or transcendental numbers? 

Only very partial results are known in the direction of Problem 3. 
A famous result of Morse and Hedlund [19] 
states that a subshift ${\cal X}$
whose complexity function satisfies $p_{\cal X}(n) \le n$ for some 
positive integer $n$ must be periodic. In that case,
it follows that ${\cal C}_{\cal X}$ is a finite set composed only
of quadratic numbers. 
Further, it is shown in [5] that for a Sturmian subshift ${\cal
X}$, that is, a subshift
with complexity $p_{\chi}(n)=n+1$ for every $n \ge 1$, the set  ${\cal
C}_{\cal X}$ 
is an uncountable set composed only by transcendental numbers. 
Theorem 4 slightly improves this result.

\proclaim Theorem 4. Let ${\cal X}$ be a subshift on $\Z_{\ge 1}$.
If the set ${\cal C}_{\cal X}$ contains a real algebraic number of degree
at least three, then the
complexity function of ${\cal X}$ satisfies 
$$
\lim_{n \to + \infty} \, p_{\cal X}(n)-n = + \infty.
$$

Linearly recurrent subshifts form
a class of particular interest of
subshifts of low complexity.
Let ${\cal X}=(X,S)$ be a subshift and $W$ be a finite
word. The cylinder associated with $W$ is, by definition, the subset 
$\langle W \rangle$ of
$X$ formed by the sequences that begin in the
word $W$. A minimal subshift 
$(X,S)$ is linearly recurrent if there exists a positive constant $c$
such that for each cylinder $\langle W \rangle$ the return time to 
$\langle W \rangle$ under $S$
is bounded by $c \vert W\vert$. 
Such dynamical systems, studied e.g. in [14], are
uniquely ergodic and have a low
complexity (in particular, they have zero
entropy), but without being necessarily trivial. 
Another contribution to Problem 3 is given by Theorem 5.

\proclaim Theorem 5. Let ${\cal X}$ be a linearly recurrent subshift
on $\Z_{\ge 1}$. Then, the set ${\cal C}_{\cal X}$ is composed only by 
quadratic or transcendental numbers. 

The proofs of Theorems 3 to 5 are postponed to Section 5.

\vskip 6mm

\centerline{\bf 4. Proofs of Theorems 1 and 2}

\vskip 6mm

The proofs of Theorems 1 and 2 rest on the following deep
result, commonly known as the Schmidt Subspace Theorem.

\proclaim Theorem A (W. {\bf M. Schmidt).} 
Let $m \ge 2$ be an integer.
Let $L_1, \ldots, L_m$ be linearly independent linear forms in
${\bf x} = (x_1, \ldots, x_m)$ with algebraic coefficients.
Let $\eps$ be a positive real number.
Then, the set of solutions ${\bf x} = (x_1, \ldots, x_m)$ 
in $\Z^m$ to the inequality
$$
\vert L_1 ({\bf x}) \ldots L_m ({\bf x}) \vert  \le
(\max\{|x_1|, \ldots , |x_m|\})^{-\eps}
$$
lies in finitely many proper subspaces of $\Q^m$.

\pro See e.g. [25] or [26].
The case $m=3$ has been established earlier
in [24]. \cqfd

Compared with the 
pioneering work [12] and the recent papers [21, 5, 13, 8],
the novelty in the present paper is that we are able to use
Theorem A with $m=4$ and not only with $m=3$, as in all of
these works.

\medskip

We further need an easy auxiliary result.

\proclaim Lemma 1. Let $\alpha = [a_0; a_1, a_2, \ldots]$ and
$\beta = [b_0; b_1, b_2, \ldots]$ be real numbers. Assume that,
for some positive integer $m$, we have $a_j = b_j$ for any $j=0, \ldots, m$.
Then, we have
$$
|\alpha - \beta| < q_m^{-2},
$$
where $q_m$ is the denominator of the convergent $[a_0; a_1, \ldots, a_m]$.

\pro Since $[a_0; a_1, \ldots, a_m] =: p_m/q_m$ is a convergent
to $\alpha$ and to $\beta$, the real numbers $\alpha - p_m/q_m$ and
$\beta - p_m/q_m$ have the same sign and are both in absolute value
less than $q_m^{-2}$, hence the lemma. \cqfd

Now, we have all the tools to establish Theorems 1 and 2.

\medskip

\noi{\bf Proof of Theorem 1.} 
Keep the notation and the hypothesis of this theorem.
Assume that the parameter $w > 1$ is fixed, as well as the 
sequence $(V_n)_{n \ge 1}$ occurring in the
definition of Condition $(*)_w$. 
Set also $s_n=\vert V_n\vert$, for any $n \ge 1$.
We want to prove that the real number
$$
\alpha:= [0; a_1, a_2, \ldots]
$$ 
is transcendental. We assume that $\alpha$ is algebraic of degree at
least three and we aim at deriving a contradiction. 
Throughout this Section, the constants implied by $\ll$ depend
only on $\alpha$.

Let $(p_{\ell}/q_{\ell})_{\ell \ge 1}$ denote the sequence of convergents 
to $\alpha$. Observe first that we have
$$
q_{\ell + 1} \ll q_{\ell}^{1.1}, \qquad
\hbox{for any $\ell \ge 1$,}   \eqno (2)
$$
by Roth's Theorem [23].

The key fact for the proof of Theorem 1 is the observation
that $\alpha$ admits infinitely many good quadratic approximants
obtained by truncating its continued fraction
expansion and completing by periodicity.
Precisely, for any positive integer $n$, we define the sequence
$(b_k^{(n)})_{k \ge 1}$ by
$$
b_{h + j s_n}^{(n)} = a_h 
\quad \hbox{for $1 \le h \le s_n$ and $j \ge 0$.} 
$$
The sequence
$(b_k^{(n)})_{k \ge 1}$ is purely periodic with period $V_n$. Set
$$
\alpha_n= [0; b_1^{(n)}, b_2^{(n)}, \ldots ]
$$
and observe that $\alpha_n$ is root of the quadratic polynomial
$$
P_n (X) := q_{s_n-1} X^2 + (q_{s_n} - p_{s_n-1}) X - p_{s_n}.
$$
By Rolle's Theorem and Lemma 1, for any positive integer $n$, 
we have
$$
|P_n (\alpha)| = |P_n (\alpha) - P_n (\alpha_n)|
\ll  q_{s_n} \, |\alpha - \alpha_n| \ll  q_{s_n} \, q_{[w s_n]}^{-2},
\eqno (3)
$$
since the first $[w s_n]$ partial quotients of $\alpha$ and $\alpha_n$
are the same. Furthermore, we clearly have
$$
|q_{s_n} \alpha - p_{s_n}| \le q_{s_n}^{-1}  \eqno (4)
$$
and we infer from (2) that
$$
|q_{s_n-1} \alpha - p_{s_n-1}| \le q_{s_n-1}^{-1} \ll q_{s_n}^{-0.9}. \eqno (5)
$$

Consider now the four linearly independent linear forms:
$$
\eqalign{
L_1(X_1, X_2, X_3, X_4) = & \alpha^2 X_2 + \alpha (X_1 - X_4) - X_3,  \cr
L_2(X_1, X_2, X_3, X_4) = & \alpha X_1 - X_3, \cr
L_3(X_1, X_2, X_3, X_4) = & X_1, \cr
L_4(X_1, X_2, X_3, X_4) = & X_2. \cr}
$$
Evaluating them on the quadruple 
$(q_{s_n}, q_{s_n-1}, p_{s_n}, p_{s_n-1})$, it follows from (3) and (4) 
that
$$
\prod_{1 \le j \le 4} \, |L_j (q_{s_n}, q_{s_n-1}, p_{s_n}, p_{s_n-1})|
\ll q_{s_n}^2 \, q_{[w s_n]}^{-2}. \eqno (6)
$$
By assumption, there exists a real number $M$ such that
$\log q_{\ell} \le \ell \, \log M$ for any positive integer $\ell$. 
Furthermore, an immediate induction shows
that $q_{\ell+2} \ge 2 \, q_{\ell}$ holds for any 
positive integer $\ell$.
Consequently, for any integer $n \ge 3$, we get
$$
{q_{[w s_n]} \over q_{s_n}} \ge \sqrt{2}^{[(w-1) s_n]-1}  
\ge q_{s_n}^{ (w - 1 - 2/s_n) (\log  \sqrt{2})/ \log M},
$$
and we infer from (6) and $w>1$ that
$$
\prod_{1 \le j \le 4} \, |L_j (q_{s_n}, q_{s_n-1}, p_{s_n}, p_{s_n-1})|
\ll q_{s_n}^{-\eps}
$$
holds for some positive real number $\eps$, when $n$ is large enough.

It then follows from
Theorem A that the points $(q_{s_n}, q_{s_n-1}, p_{s_n}, p_{s_n-1})$
lie in a finite number of proper subspaces of $\Q^4$. 
Thus, there exist a non-zero integer quadruple $(x_1,x_2,x_3,x_4)$ and
an infinite set of distinct positive integers ${\cal N}_1$ such that
$$
x_1 q_{s_n} + x_2 q_{s_n - 1} + x_3 p_{s_n} + x_4 p_{s_n - 1} = 0,  \eqno (7)
$$
for any $n$ in ${\cal N}_1$.
Observe that $(x_2, x_4) \not= (0, 0)$, since, otherwise, 
by letting $n$
tend to infinity along ${\cal N}_1$ in (7), 
we would get that the real number
$\alpha$ is rational. Dividing (7) by $q_{s_n}$, we obtain
$$
x_1 + x_2 {q_{s_n - 1} \over q_{s_n}} + x_3 {p_{s_n} \over q_{s_n}} 
+ x_4 { p_{s_n - 1} \over q_{s_n - 1} } \cdot {q_{s_n - 1} \over q_{s_n}}
= 0.  \eqno (8)
$$
By letting $n$ tend to infinity along ${\cal N}_1$ in (8),
we get that
$$
\beta := \lim_{{\cal N}_1 \ni n \to + \infty} \, {q_{s_n - 1} \over q_{s_n}}
= - {x_1 + x_3 \alpha \over x_2 + x_4 \alpha}.
$$
Furthermore, observe that, for any $n$ in ${\cal N}_1$, we have
$$
\biggl| \beta - {q_{s_n - 1} \over q_{s_n}} \biggr| =
\biggl| {x_1 + x_3 \alpha \over x_2 + x_4 \alpha} - {x_1 + x_3 p_{s_n}/q_{s_n}
\over x_2 + x_4 p_{s_n - 1} / q_{s_n - 1}} \biggr| 
\ll {1 \over q_{s_n-1}^2} \ll {1 \over q_{s_n}^{1.8}},
\eqno (9)
$$
by (4) and (5).
Since $q_{s_n-1}$ and $q_{s_n}$ are coprime and
$s_n$ tends to infinity when $n$ tends to infinity along ${\cal N}_1$,
this implies that $\beta$ is irrational.

Consider now the three linearly independent
linear forms:
$$
L'_1(Y_1, Y_2, Y_3) = \beta Y_1  - Y_2, \quad 
L'_2(Y_1, Y_2, Y_3) = \alpha Y_1 - Y_3, \quad 
L'_3(Y_1, Y_2, Y_3) = Y_1. 
$$
Evaluating them on the triple 
$(q_{s_n}, q_{s_n-1}, p_{s_n})$ with $n \in {\cal N}_1$, 
we infer from (4) and (9) that
$$
\prod_{1 \le j \le 3} \, |L'_j (q_{s_n}, q_{s_n-1}, p_{s_n})|
\ll  q_{s_n}^{-0.8}.
$$
It then follows from
Theorem A that the points $(q_{s_n}, q_{s_n-1}, p_{s_n})$ 
with $n \in {\cal N}_1$ lie in a finite number of proper subspaces of $\Q^3$. 
Thus, there exist a non-zero integer triple $(y_1, y_2, y_3)$ and
an infinite set of distinct positive integers ${\cal N}_2$ such that
$$
y_1 q_{s_n} + y_2 q_{s_n - 1} + y_3 p_{s_n}  = 0,  \eqno (10)
$$
for any $n$ in ${\cal N}_2$.
Dividing (10) by $q_{s_n}$ and letting $n$ tend to
infinity along ${\cal N}_2$, we get
$$
y_1 + y_2 \beta + y_3 \alpha = 0.  \eqno (11)
$$

To obtain another equation linking $\alpha$ and $\beta$, we
consider the three linearly independent
linear forms:
$$
L''_1(Z_1, Z_2, Z_3) = \beta Z_1 - Z_2, \quad 
L''_2(Z_1, Z_2, Z_3) = \alpha Z_2 - Z_3, \quad 
L''_3(Z_1, Z_2, Z_3) = Z_1. 
$$
Evaluating them on the triple 
$(q_{s_n}, q_{s_n-1}, p_{s_n-1})$ with $n \in {\cal N}_1$, 
we infer from (5) and (9) that
$$
\prod_{1 \le j \le 3} \, |L''_j (q_{s_n}, q_{s_n-1}, p_{s_n-1})|
\ll  q_{s_n}^{-0.7}.
$$
It then follows from
Theorem A that the points $(q_{s_n}, q_{s_n-1}, p_{s_n-1})$ 
with $n \in {\cal N}_1$ lie in a finite number of proper subspaces of $\Q^3$. 
Thus, there exist a non-zero integer triple $(z_1, z_2, z_3)$ and
an infinite set of distinct positive integers ${\cal N}_3$ such that
$$
z_1 q_{s_n} + z_2 q_{s_n - 1} + z_3 p_{s_n-1}  = 0,  \eqno (12)
$$
for any $n$ in ${\cal N}_3$.
Dividing (12) by $q_{s_n-1}$ and letting $n $ tend to
infinity along ${\cal N}_3$, we get
$$
{z_1 \over \beta} + z_2   + z_3 \alpha = 0.  \eqno (13)
$$
Observe that $y_2\not= 0$ since $\alpha$ is irrational.
We infer from (11) and (13) that
$$
(z_3 \alpha + z_2) (y_3 \alpha + y_1) = y_2 z_1.
$$
If $y_3 z_3 = 0$, then (11) and (13) yield that $\beta$ is rational,
which is a contradiction. Consequently, $y_3 z_3 \not= 0$ and
$\alpha$ is a quadratic real number, which is again a contradiction. 
This completes the proof of the second assertion of the theorem. 

\medskip

It then remains for us to explain why we can drop the assumption on the
sequence $(q_{\ell}^{1/\ell})_{\ell \ge 1}$ when $w$ is sufficiently large.
We return to the beginning of the proof, and we
assume that $w \ge 2$. Using well-known facts from the theory
of continuants (see e.g. [20]), inequality (3) becomes
$$
|P_n (\alpha)| \ll q_{s_n} \, q_{2 s_n}^{-2} \ll q_{s_n} \, q_{s_n}^{-4}
\ll q_{s_n}^{-3} \ll {\rm H} (P_n)^{-3},
$$
where ${\rm H} (P_n)$ denotes the height of the polynomial $P_n$, that
is, the maximum of the absolute values of its coefficients.
By the main result from [24] (or by using Theorem A with $m=3$
and the linear forms $\alpha^2 X_2 + \alpha X_1 + X_0$, $X_2$ and $X_1$), 
this immediately implies that
$\alpha$ is transcendental. \cqfd

\bigskip \bigskip

\noi{\bf Proof of Theorem 2.} 
Assume that the parameters $w$ and $w'$ are fixed, as well as the 
sequences $(U_n)_{n \ge 1}$ and $(V_n)_{n \ge 1}$ occurring in the
definition of Condition $(**)_{w, w'}$. 
Without any loss of generality, we add in the statement of
Condition $(**)_{w, w'}$ the following two assumptions:
\medskip

\item{\rm (iv)} The sequence $(\vert U_n\vert)_{n \ge 1}$ is unbounded;

\smallskip

\item{\rm (v)} For any $n \ge 1$, the last letter of the word $U_n$
differs from the last letter of the word $V_n$.

\medskip

We point out that the conditions $(iv)$ and $(v)$ 
do not at all restrict the generality.
Indeed, if $(iv)$ is not fulfilled by a sequence ${\bf a}$
satisfying $(i)-(iii)$ of Condition $(**)_{w, w'}$, then 
the desired result follows from Theorem 1. 
To see that $(v)$ does not cause any trouble,
we make the following observation. Let $a$ be a letter
and $U$ and $V$ be two words such that ${\bf a}$ begins with
$Ua(Va)^w$. Then, ${\bf a}$ also begins with $U (aV)^w$ and 
we have trivially $|U|/|aV| \le |Ua|/|Va|$.

Set $r_n=\vert U_n\vert$ and $s_n=\vert V_n\vert$, for any $n \ge 1$.
We want to prove that the real number
$$
\alpha:= [0; a_1, a_2, \ldots]
$$ 
is transcendental. We assume that $\alpha$ is algebraic of degree at
least three and we aim at deriving a contradiction. 
Let $(p_{\ell}/q_{\ell})_{\ell \ge 1}$ denote the sequence of convergents 
to $\alpha$.

Let $n$ be a positive integer. Since $w>1$ and $r_n \le w' s_n$, we get
$$
{2 r_n + s_n \over r_n + w s_n} \le {2 w' s_n + s_n \over w' s_n + w s_n} 
= {2 w' + 1 \over w' + w} < {\log m \over \log M},
$$
by (1). Consequently, there
exist positive real numbers $\eta$ and $\eta'$ with $\eta < 1$ such that
$$
(1 + \eta) (1 + \eta') (2 r_n + s_n) \log M <
(1 - \eta') (r_n + w s_n) \log m,  \eqno (14)
$$
for any $n \ge 1$. Notice that we have
$$
q_{\ell + 1} \ll q_{\ell}^{1 + \eta}, \qquad
\hbox{for any $\ell \ge 1$,}   \eqno (15)
$$
by Roth's Theorem [23].

As for the proof of Theorem 1, we observe
that $\alpha$ admits infinitely many good quadratic approximants
obtained by truncating its continued fraction
expansion and completing by periodicity.
Precisely, for any positive integer $n$, we define the sequence
$(b_k^{(n)})_{k \ge 1}$ by
$$
\eqalign{
b_h^{(n)} & = a_h \quad \hbox{for $1 \le h \le r_n + s_n$,} \cr
b_{r_n + h + j s_n}^{(n)} & = a_{r_n + h} 
\quad \hbox{for $1 \le k \le s_n$ and $j \ge 0$.} \cr}
$$
The sequence
$(b_k^{(n)})_{k \ge 1}$ is eventually periodic, with preperiod $U_n$
and with period $V_n$. Set
$$
\alpha_n= [0; b_1^{(n)}, b_2^{(n)}, \ldots ]
$$
and observe that $\alpha_n$ is root of the quadratic polynomial
$$
\eqalign{
P_n (X) & := (q_{r_n-1} q_{r_n+s_n} - q_{r_n} q_{r_n+s_n-1} ) X^2  \cr
& - (q_{r_n-1} p_{r_n+s_n} - q_{r_n} p_{r_n+s_n-1} 
+ p_{r_n-1} q_{r_n+s_n} - p_{r_n} q_{r_n+s_n-1}) X \cr
& + (p_{r_n-1} p_{r_n+s_n} - p_{r_n} p_{r_n+s_n-1}). \cr}
$$
For any positive integer $n$, we infer from Rolle's Theorem 
and Lemma 1 that
$$
|P_n (\alpha)| = |P_n (\alpha) - P_n (\alpha_n)|
\ll \, q_{r_n} \,  q_{r_n + s_n} \, |\alpha - \alpha_n| \ll 
q_{r_n} \, q_{r_n + s_n} \, q_{r_n + [w s_n]}^{-2},
\eqno (16)
$$
since the first $r_n + [w s_n]$ partial quotients of $\alpha$ and $\alpha_n$
are the same. Furthermore, by (15), we have
$$
|(q_{r_n-1} q_{r_n+s_n} - q_{r_n} q_{r_n+s_n-1} ) \alpha 
- (q_{r_n-1} p_{r_n+s_n} - q_{r_n} p_{r_n+s_n-1})| \ll
q_{r_n} \, q_{r_n + s_n}^{-1 + \eta}  \eqno (17)
$$
and
$$
|(q_{r_n-1} q_{r_n+s_n} - q_{r_n} q_{r_n+s_n-1} ) \alpha 
- (p_{r_n-1} q_{r_n+s_n} - p_{r_n} q_{r_n+s_n-1})| \ll
q_{r_n}^{-1 + \eta} \, q_{r_n + s_n}.  \eqno (18)
$$
We have as well the obvious upper bound
$$
|q_{r_n-1} q_{r_n+s_n} - q_{r_n} q_{r_n+s_n-1} | 
\le  q_{r_n} \, q_{r_n + s_n}.
\eqno (19)
$$

Consider now the four linearly independent linear forms:
$$
\eqalign{
L_1(X_1, X_2, X_3, X_4) = & \alpha^2 X_1 - \alpha (X_2 + X_3) + X_4,  \cr
L_2(X_1, X_2, X_3, X_4) =& \alpha X_1 - X_2, \cr
L_3(X_1, X_2, X_3, X_4) =& \alpha X_1 - X_3, \cr
L_4(X_1, X_2, X_3, X_4) =& X_1. \cr}
$$
Evaluating them on the quadruple 
$$
\eqalign{
{\underline z_n} := (q_{r_n-1} q_{r_n+s_n} & - q_{r_n} q_{r_n+s_n-1}, 
q_{r_n-1} p_{r_n+s_n} - q_{r_n} p_{r_n+s_n-1},  \cr
& p_{r_n-1} q_{r_n+s_n} - p_{r_n} q_{r_n+s_n-1}, 
p_{r_n-1} p_{r_n+s_n} - p_{r_n} p_{r_n+s_n-1}), \cr}
$$ 
it follows from (16), (17), (18), and (19) that
$$
\prod_{1 \le j \le 4} \, |L_j ({\underline z_n})|
\ll q_{r_n}^{2 + \eta} \, q_{r_n + s_n}^{2 + \eta} \, q_{r_n + [w s_n]}^{-2} 
\ll ( q_{r_n} \, q_{r_n + s_n})^{-\eta}
(q_{r_n}^{1 + \eta} \, q_{r_n + s_n}^{1 + \eta} \, q_{r_n + [w s_n]}^{-1})^2.
$$
Assuming $n$ sufficiently large, we have
$$
q_{r_n} \le M^{(1 + \eta') r_n}, \qquad
q_{r_n + s_n} \le M^{(1 + \eta') (r_n+s_n)}, \quad {\rm and}
\quad q_{r_n + [w s_n]} \ge m^{(1 - \eta') (r_n + w s_n)},
$$
with $\eta'$ as in (14). Consequently, we get
$$
(q_{r_n}^{1 + \eta} \, q_{r_n + s_n}^{1 + \eta} \, q_{r_n + [w s_n]}^{-1})
\le M^{(1+\eta)(1 + \eta')(2 r_n + s_n)} \, m^{- (1 - \eta') (r_n + w s_n)}
\le 1,
$$
by (14). Thus, we get the upper bound
$$
\prod_{1 \le j \le 4} \, |L_j ({\underline z_n})|
\ll ( q_{r_n} \, q_{r_n + s_n})^{-\eta}
$$
for any positive integer $n$.

It then follows from
Theorem A that the points ${\underline z_n}$
lie in a finite number of proper subspaces of $\Q^4$. 
Thus, there exist a non-zero integer quadruple $(x_1,x_2,x_3,x_4)$ and
an infinite set of distinct positive integers ${\cal N}_1$ such that
$$
\eqalign{
& x_1 (q_{r_n-1} q_{r_n+s_n} - q_{r_n} q_{r_n+s_n-1}) 
+ x_2 (q_{r_n-1} p_{r_n+s_n} - q_{r_n} p_{r_n+s_n-1})  \cr
& + x_3 (p_{r_n-1} q_{r_n+s_n} - p_{r_n} q_{r_n+s_n-1}) 
+ x_4 (p_{r_n-1} p_{r_n+s_n} - p_{r_n} p_{r_n+s_n-1}) = 0, \cr} \eqno (20)
$$
for any $n$ in ${\cal N}_1$.

Divide (20) by $q_{r_n} \, q_{r_n+s_n -1}$ and observe that
$p_{r_n}/q_{r_n}$ and $p_{r_n+s_n}/q_{r_n+s_n}$ tend to $\alpha$
as $n$ tends to infinity along ${\cal N}_1$. Taking the limit, we get that
either
$$
x_1 + (x_2 + x_3) \alpha + x_4 \alpha^2 = 0 \eqno (21)
$$
or
$$
{q_{r_n - 1} q_{r_n+s_n}  \over q_{r_n} q_{r_n+s_n - 1}  }
\qquad 
\hbox{tends to $1$ as $n$ tends to infinity along ${\cal N}_1$} \eqno (22)
$$
must hold. In the former case, since $\alpha$ is irrational
and not quadratic, we get that $x_1 = x_4 = 0$ and $x_2 = - x_3$. 
Then, $x_2$ is non-zero and, for any $n$ in ${\cal N}_1$, we have
$q_{r_n-1} p_{r_n+s_n} - q_{r_n} p_{r_n+s_n-1} =
p_{r_n-1} q_{r_n+s_n} - p_{r_n} q_{r_n+s_n-1}$.
Thus, the polynomial
$P_n(X)$ can simply be expressed as
$$
\eqalign{
P_n (X) & := (q_{r_n-1} q_{r_n+s_n} - q_{r_n} q_{r_n+s_n-1} ) X^2  \cr
& - 2 (q_{r_n-1} p_{r_n+s_n} - q_{r_n} p_{r_n+s_n-1}) X 
+ (p_{r_n-1} p_{r_n+s_n} - p_{r_n} p_{r_n+s_n-1}). \cr}
$$
Consider now the three linearly independent linear forms:
$$
\eqalign{
L'_1 (Y_1, Y_2, Y_3) = & \alpha^2 X_1 - 2 \alpha X_2 + X_3,  \cr
L'_2 (Y_1, Y_2, Y_3) =  & \alpha  X_1 - X_2, \cr
L'_3 (Y_1, Y_2, Y_3) =  & X_1. \cr}
$$
Evaluating them on the triple 
$$
\eqalign{
{\underline z'_n} := (q_{r_n-1} q_{r_n+s_n} & - q_{r_n} q_{r_n+s_n-1}, 
q_{r_n-1} p_{r_n+s_n} - q_{r_n} p_{r_n+s_n-1},  \cr
& p_{r_n-1} p_{r_n+s_n} - p_{r_n} p_{r_n+s_n-1}), \cr}
$$ 
it follows from (16), (17) and (19) that
$$
\prod_{1 \le j \le 3} \, |L'_j ({\underline z'_n})|
\ll q_{r_n}^3 \, q_{r_n + s_n}^{1 + \eta} \, q_{r_n + [w s_n]}^{-2} 
\ll q_{r_n}^2 \, q_{r_n + s_n}^{2 + \eta} \, q_{r_n + [w s_n]}^{-2} 
\ll ( q_{r_n} \, q_{r_n + s_n})^{-\eta},
$$
by the above computation.

It then follows from
Theorem A that the points ${\underline z'_n}$
lie in a finite number of proper subspaces of $\Q^3$. 
Thus, there exist a non-zero integer triple $(x'_1,x'_2,x'_3)$ and
an infinite set of distinct positive integers ${\cal N}_2$
included in ${\cal N}_1$ such that
$$
\eqalign{
x'_1 (q_{r_n-1} q_{r_n+s_n} & - q_{r_n} q_{r_n+s_n-1}) 
+ x'_2 (q_{r_n-1} p_{r_n+s_n} - q_{r_n} p_{r_n+s_n-1})  \cr
& + x'_3 (p_{r_n-1} p_{r_n+s_n} - p_{r_n} p_{r_n+s_n-1}) = 0, \cr} 
\eqno (23)
$$
for any $n$ in ${\cal N}_2$. 

Divide (23) by $q_{r_n} \, q_{r_n+s_n -1}$ and observe that
$p_{r_n}/q_{r_n}$ and $p_{r_n+s_n}/q_{r_n+s_n}$ tend to $\alpha$
as $n$ tends to infinity along ${\cal N}_2$. Taking the limit, we get that
either
$$
x'_1 + x'_2 \alpha + x'_3 \alpha^2 = 0  \eqno (24)
$$
or
$$
{q_{r_n - 1} q_{r_n+s_n}  \over q_{r_n} q_{r_n+s_n - 1}  }
\qquad 
\hbox{tends to $1$ as $n$ tends to infinity along ${\cal N}_2$} \eqno (25)
$$
must hold. In the former case, we have a contradiction
since $\alpha$ is irrational and not quadratic.

Consequently, to conclude the proof of our theorem, 
it is enough to derive a contradiction
from (22) (resp. from (25)), assuming that (21) (resp. (24)) does not
hold. To this end, we observe that (20) (resp. (23))
allows us to control the speed of
convergence of $Q_n := ( q_{r_n - 1} q_{r_n+s_n} )/  (q_{r_n} q_{r_n+s_n - 1})$
to $1$ along ${\cal N}_1$ (resp. along ${\cal N}_2$). 

Thus, we assume that the quadruple $(x_1, x_2, x_3, x_4)$ 
obtained after the first
application of Theorem A satisfies 
$x_1 + (x_2 + x_3) \alpha + x_4 \alpha^2 \not= 0$.
Dividing (20) by $q_{r_n} q_{r_n+s_n - 1}$, we get
$$
\eqalign{
& x_1 (Q_n - 1) 
+ x_2 \biggl(Q_n { p_{r_n+s_n} \over q_{r_n+s_n}}
- {p_{r_n+s_n-1} \over q_{r_n+s_n-1} } \biggr) 
+ x_3 \biggl(Q_n { p_{r_n-1} \over q_{r_n-1}}  
- {p_{r_n} \over q_{r_n} } \biggr) \cr
& + x_4 \biggl( Q_n {p_{r_n-1} \over q_{r_n-1}}
{p_{r_n+s_n} \over q_{r_n + s_n}} - {p_{r_n} \over q_{r_n}}
{ p_{r_n+s_n-1} \over q_{r_n + s_n - 1}} \biggr) = 0, \cr} \eqno (26)
$$
for any $n$ in ${\cal N}_1$. To shorten the notation, for any $\ell \ge 1$,
we put $R_\ell := \alpha - p_\ell/q_\ell$. We rewrite (26) as
$$
\eqalign{
& x_1 (Q_n - 1) 
+ x_2 \bigl(Q_n (\alpha - R_{r_n+s_n})
- (\alpha - R_{r_n+s_n-1}) \bigr) 
+ x_3 \bigl(Q_n (\alpha - R_{r_n-1})  - (\alpha - R_{r_n}) \bigr) \cr
& + x_4 \bigl( Q_n (\alpha - R_{r_n-1})
(\alpha - R_{r_n+s_n}) - (\alpha - R_{r_n})
(\alpha -R_{r_n+s_n-1} ) \bigr) = 0. \cr} 
$$
This yields
$$
\eqalign{
& (Q_n-1)  \bigl( x_1 + (x_2 + x_3) \alpha + x_4 \alpha^2 \bigr) \cr
= x_2 Q_n R_{r_n+s_n} - x_2 R_{r_n+s_n-1} + x_3 Q_n & R_{r_n-1} - x_3 R_{r_n} 
- x_4 Q_n R_{r_n-1} R_{r_n+s_n} \cr 
+ x_4 R_{r_n} R_{r_n+s_n-1} +  \alpha (x_4 Q_n R_{r_n-1} 
+ x_4 Q_n & R_{r_n+s_n} - x_4  R_{r_n} -
x_4 R_{r_n+s_n-1}). \cr} \eqno (27)
$$
Observe that $|R_\ell| \le 1/q_\ell^2$ for any $\ell \ge 1$. Furthermore, for 
$n$ large enough, we have $1/2 \le Q_n \le 2$, by our assumption (22).
Consequently, we derive from (27) that
$$
|(Q_n-1) (x_1 + (x_2 + x_3) \alpha + x_4 \alpha^2)| \ll
|R_{r_n-1}| \ll q_{r_n-1}^{-2}.
$$
Since we have assumed that (21) does not hold, we get
$$
|Q_n - 1| \ll q_{r_n-1}^{-2}. \eqno (28)
$$

On the other hand, observe that the rational number
$Q_n$ is the quotient of the two continued fractions
$[a_{r_n+s_n}; a_{r_n+s_n-1}, \ldots , a_1]$ and
$[a_{r_n}; a_{r_n-1}, \ldots , a_1]$. By assumption (v) from
Condition $(**)_{w, w'}$, we have
$a_{r_n+s_n} \not= a_{r_n}$, thus either $a_{r_n+s_n}-a_{r_n}\geq 1$ or
$a_{r_n}-a_{r_n+s_n}\geq 1$ holds. A simple calculation then shows
that
$$
|Q_n - 1| \gg a_{r_n}^{-1} \, \min \{ a_{r_n + s_n - 1}^{-1} +
a_{r_n - 2}^{-1}, a_{r_n + s_n - 2}^{-1} + a_{r_n - 1}^{-1} \} \gg
a_{r_n}^{-1} \, q_{r_n - 1}^{-1}, 
$$
since $q_{r_n - 1} \ge \max\{a_{r_n - 1}, a_{r_n - 2}\}$. Combined
with (28), this gives $a_{r_n} \gg q_{r_n - 1}$ and
$$
q_{r_n} \ge a_{r_n} q_{r_n-1} \gg q_{r_n-1}^2. \eqno (29)
$$
Since $\eta < 1$ and (29) holds for infinitely many $n$, 
we get a contradiction with (15).

We derive a contradiction from (25) in an entirely similar way.
This completes the proof of our theorem. \cqfd

\vskip 6mm

\centerline{\bf 5. Proofs of Theorem 3 to 5}

\vskip 6mm

Before establishing Theorems 3 to 5,
we state an easy, but useful, auxiliary result. 

\proclaim Lemma 2. 
Let $\sigma$ be a prolongable morphism defined on a finite alphabet ${\cal A}$. 
Let ${\bf a}$ be the associated fixed point 
and $a$ be the first letter of ${\bf a}$. 
Then, there exists a positive constant $c$ such that, for any
positive integer $n$ and any letter $b$ occurring
in ${\bf a}$, we have
$\vert\sigma^n(a)\vert\geq c\vert\sigma^n(b)\vert$.

\noi {\bf Proof of Lemma 2}. 
Without loss of generality,
we may assume that ${\cal A}$ is exactly the set
of letters occurring in ${\bf a}$. 
Let $b$ be in ${\cal A}$.
Since ${\bf a}$ is obtained as 
the limit $\lim_{n\to + \infty}\sigma^n(a)$, 
there exists an integer $n_b$ such
that the word $\sigma^{n_b}(a)$ contains the letter $b$. Set 
$$
s=\max_{b\in{\cal A}}\left\{\vert\sigma(b)\vert\right\}
\;\;\hbox{ and }\;\;n_0=\max_{b\in{\cal A}}\left\{n_b\right\}.
$$ 
Let $n$ be a positive integer.
If $n\leq n_0$, then we have 
$$
\vert\sigma^n(a)\vert\geq 1\geq {1\over{s^{n_0}}} \,
\vert\sigma^n(b)\vert.
$$ 
If $n>n_0$, then we get
$$
\vert\sigma^n(a)\vert=\vert\sigma^{n-n_0}(\sigma^{n_0}(a))
\vert\geq\vert\sigma^{n-n_0}(b)
\vert\geq {1\over{s^{n_0}}} \, \vert\sigma^n(b)\vert,
$$ 
and the lemma follows by taking $c=s^{-n_0}$.
\cqfd

\noi {\bf Proof of Theorem 3}. 
Let us assume that ${\bf b}$ is a sequence generated by a recurrent morphism 
$\sigma$ and that ${\bf b}$ is not eventually periodic. 
There exists a fixed point ${\bf a}$  of $\sigma$ and a coding $\varphi$ 
such that ${\bf b}=\varphi({\bf a})$. 
By assumption, the first letter $a$
which occurs in ${\bf a}$ should appear at least twice. Thus, there
exists a finite (possibly empty) word $W$ such that $aWa$ is a
prefix of the word ${\bf a}$.  
We check that the assumptions of Theorem 1 are satisfied
by ${\bf b}$ with the sequence $(V_n)_{n \ge 1}$ 
defined by $V_n = \varphi(\sigma^{n} (aW))$  
for any $n \ge 1$. Indeed, by Lemma 2, there exists a positive
rational number $c$, depending only on $\sigma$ and $W$, such that 
$$
\vert\sigma^n(a)\vert\geq  c \, \vert\sigma^n(aW)\vert. 
$$
This implies that $\varphi(\sigma^n(aWa))$ begins in  
$(\varphi(\sigma^n(aW)))^{1+c}$. 
Since $\varphi(\sigma^n(aWa))$ is a prefix of ${\bf a}$,  
we get that the sequence ${\bf b}$ satisfies Condition 
$(*)_{1+c}$. We conclude by applying Theorem 1.
\cqfd

\noi {\bf Proof of Theorem 4}. Let ${\cal X} = (X, S)$ 
be a subshift such that
$p_{\cal X}(n)-n$ is bounded and let 
$\alpha$ be an element of ${\cal C}_{\cal X}$. 
By definition of the set  ${\cal C}_{\cal X}$, there exists a
sequence ${\bf a}=(a_\ell)_{\ell \geq 1}$ in $X$ such that
$\alpha=[0;a_1,a_2,\ldots]$.
 
First, assume that the complexity function of the sequence
${\bf a}$ satisfies $p_{\bf a} (n)\leq n$ for some $n$. 
It follows from a theorem of Morse and Hedlund [19] 
that ${\bf a}$ is eventually periodic, thus $\alpha$ is a quadratic
number. 

Now, assume that $p_{\cal X} (n) > n$ for every integer $n$. 
Since $p_{\cal X} (n) - n$ is bounded,
there exist two positive
integers $n_0$ and $a$ such that  $p_{\cal X} (n)=n+a$ for all $n\geq n_0$
(see for instance [4]). 
This implies (see e.g. [10])
that there exist a finite word $W$, a non-erasing morphism
$\phi$ and a Sturmian sequence ${\bf u}$ such that ${\bf a}=W\phi({\bf
u})$. Since ${\bf u}$ begins in arbitrarily long
squares (this is proved in [5])
and since $\phi$ is a non-erasing morphism, it follows that 
$\phi({\bf u})$ also begins in arbitrarily long squares, hence, it
satisfies Condition
$(*)_2$. We then infer from Theorem 1 that the real number
$\alpha'=[0;a_{\vert W\vert+1},a_{\vert W\vert+2},a_{\vert
W\vert+3},\ldots]$ is transcendental. It immediately follows that 
$\alpha$ is a transcendental number, concluding the proof
of the theorem. \cqfd

\noi {\bf Proof of Theorem 5}. 
Let ${\cal X}=(X,S)$ be a linearly recurrent subshift and let $\alpha$ 
be an element of ${\cal C}_{\cal X}$. By the definition of the set  
${\cal C}_{\cal X}$, there exists a
sequence ${\bf a}=(a_\ell)_{\ell \geq 1}$ in $X$ such that
$\alpha=[0;a_1,a_2,\ldots]$. 
By assumption, there exists a positive integer $k$ such that the gap between
two consecutive occurrences 
in ${\bf a}$ of any factor $W$ of length $n$ is at most
$kn$. For every positive integer $n$, let $U_n$ denote the prefix of
length $n$ of ${\bf a}$ and let $W_n$ be the word of length $kn$
defined by $U_{(k+1)n}=U_n W_n$. Since, by assumption, $U_n$ has at least one
occurrence in the word $W_n$, there exist two (possibly empty)
finite words $A_n$ and $B_n$ such that $W_n=A_n U_n B_n$. It follows
that $U_n A_n U_n$ is a prefix of ${\bf a}$ and, moreover,
$U_n A_n U_n=(U_n A_n)^w$ for some rational number $w$
with $w\geq 1+ 1/k$. Then, either ${\bf
a}$ is eventually periodic (in which case $\alpha$ is a quadratic
number) or ${\bf a}$ satisfies the Conditon $(*)_{1+1/k}$ and the
transcendence of $\alpha$ follows from Theorem 1, concluding the proof.
\cqfd

\vskip 4mm

\noindent {\bf Acknowledgements.} We would like to warmly thank  
Jean-Paul Allouche for many useful remarks. The first author is also most grateful to 
Val\'erie Berth\'e for her constant help and support.


\vskip 8mm

\centerline{\bf References}

\vskip 6mm

\item{[1]}
B. Adamczewski \& Y. Bugeaud,
{\it On the complexity of algebraic numbers I. Expansions
in integer bases}. Preprint.

\item{[2]}
B. Adamczewski, Y. Bugeaud \& J. L. Davison,
{\it Transcendental continued fractions}. In preparation.

\item{[3]}
B. Adamczewski, Y. Bugeaud \&  F. Luca,
{\it Sur la complexit\'e des nombres alg\'ebriques},
C. R. Acad. Sci. Paris 339 (2004), 11--14.

\item{[4]}
J.-P. Allouche,
{\it Nouveaux r\'esultats de transcendance de r\'eels \`a 
d\'eveloppements non al\'eatoire},
Gaz. Math. 84 (2000), 19--34.

\item{[5]}
J.-P. Allouche, J. L. Davison, M. Queff\'elec \& L. Q. Zamboni,
{\it Transcendence of Sturmian or morphic continued fractions},
J. Number Theory 91 (2001), 39--66.

\item{[6]}
J.-P. Allouche \& J. O. Shallit, 
Automatic Sequences: Theory, Applications, Generalizations, 
Cambridge University Press, Cambridge, 2003.

\item{[7]}
D. H. Bailey, J. M. Borwein, R. E. Crandall \& C. Pomerance,
{\it On the binary expansions of algebraic numbers},
J. Th\'eor. Nombres Bordeaux. To appear.

\item{[8]}
C. Baxa,
{\it Extremal values of continuants and transcendence of certain
continued fractions},
Adv. in Appl. Math. 32 (2004), 754--790.

\item{[9]}
   \'E. Borel,
 {\it Sur les chiffres d\'ecimaux de $\sqrt{2}$ et divers
              probl\`emes de probabilit\'es en cha\^\i ne},
 C.~ R.~ Acad.~ Sci.~ Paris 230 (1950), 591--593.

\item{[10]}
J. Cassaigne, 
{\it Sequences with grouped factors}. In: DLT'97, Developments in
Language Theory III, Thessaloniki, Aristotle University of
Thessaloniki (1998), 211--222.

\item{[11]}
 A. Cobham,
 {\it Uniform tag sequences},
Math. Systems Theory 6 (1972), 164--192.

\item{[12]}
J. L. Davison,
{\it A class of transcendental numbers with bounded partial quotients}.
In R. A. Mollin, ed., Number Theory and Applications, pp. 365--371, Kluwer
Academic Publishers, 1989.

\item{[13]}
J. L. Davison,
{\it Continued fractions with bounded partial quotients},
Proc. Edinburgh Math. Soc. 45 (2002), 653--671.

\item{[14]}
F. Durand,
{\it Linearly recurrent subshifts have a finite number of
              non-periodic subshift factors}, 
Erg. Th. Dyn. Syst. 20  (2000), 1061--1078. Corrigendum and addendum
23 (2003), 663--669.

\item{[15]}
 J. Hartmanis \& R. E. Stearns,
{\it On the computational complexity of algorithms},
Trans. Amer. Math. Soc. 117 (1965), 285--306.

\item{[16]}
A. Ya. Khintchine,
Continued fractions, Gosudarstv. Izdat. Tehn.-Theor. Lit. 
Moscow-Leningrad, 2nd edition, 1949 (in Russian).

\item{[17]}
S. Lang,
Introduction to Diophantine Approximations, Sprin\-ger-Verlag (1995).

\item{[18]}
M. Mend\`es France,
{\it Principe de la sym\'etrie perturb\'ee}. 
In: {S\'eminaire de Th\'eorie des
Nombres, Paris 1979-80}, M.-J. Bertin (\'ed.), Birkh\"auser, Boston, 1981,
pp. 77--98.

\item{[19]}
M. Morse \& G. A. Hedlund,
{\it Symbolic dynamics},
Amer. J. Math. 60 (1938), 815--866.

\item{[20]}
O. Perron,
Die Lehre von den Ketterbr\"uchen.
Teubner, Leipzig, 1929.

\item{[21]}
M. Queff\'elec,
{\it Transcendance des fractions continues de Thue--Morse},
J. Number Theory 73 (1998), 201--211.

\item{[22]}
M. Queff\'elec,
{\it Irrational number with automaton-generated continued fraction expansion}.
In: J.-M. Gambaudo, P. Hubert, P. Tisseur, and S. Vaienti, editors, 
{\it Dynamical Systems: From Crystal to Chaos}, World Scientific,
(2000), 190--198.

\item{[23]}
K. F. Roth,
{\it Rational approximations to algebraic numbers},
Mathematika 2 (1955), 1--20. Corrigendum {\it ibid}, 168.

\item{[24]}
 W. M. Schmidt,
{\it On simultaneous approximations of two algebraic numbers by rationals},
Acta Math. 119 (1967), 27--50.

\item{[25]}
 W. M. Schmidt,
{\it Norm form equations},
Ann. of Math. 96 (1972), 526--551.
  
\item{[26]}
W. M. Schmidt, 
{\it Diophantine approximation},
Lecture Notes in Mathematics 785, Springer, Berlin, 1980.

\item{[27]}
J. O. Shallit, 
{\it Real numbers with bounded partial quotients}, 
Enseign. Math. 38 (1992), 151--187.

\item{[28]}
A. M. Turing,
{\it On computable numbers, with an application to the 
Entscheidungsproblem},
Proc. London Math. Soc. 42 (1937), 230--265.

\item{[29]}
M. Waldschmidt,
{\it Un demi-si\`ecle de transcendance}. In:
Development of mathematics 1950--2000, pp. 1121--1186, 
Birkh\"auser, Basel (2000).

\bigskip\bigskip

\noindent Boris Adamczewski   \hfill{Yann Bugeaud}

\noindent   CNRS, Institut Camille Jordan  
\hfill{Universit\'e Louis Pasteur}

\noindent   Universit\'e Claude Bernard Lyon 1 
\hfill{U. F. R. de math\'ematiques}

\noindent   B\^at. Braconnier, 21 avenue Claude Bernard
 \hfill{7, rue Ren\'e Descartes}

\noindent   69622 VILLEURBANNE Cedex (FRANCE)   
\hfill{67084 STRASBOURG Cedex (FRANCE)}

\vskip2mm
 
\noindent {\tt Boris.Adamczewski@math.univ-lyon1.fr}
\hfill{{\tt bugeaud@math.u-strasbg.fr}}

\bye